\documentclass[12pt]{amsart}
\usepackage{amssymb}
\usepackage[dvips]{graphics}

\ifx\mathrm\undefined\let\mathrm\rm\fi
\ifx\mathbf\undefined\let\mathbf\bf\fi
\ifx\mathfrak\undefined\let\mathfrak\frak\fi
\ifx\mathcal\undefined\let\mathcal\cal\fi
\ifx\mathbb\undefined\let\mathbb\Bbb\fi
\ifx\emph\undefined\let\emph\it\fi
\font\bb=msbm10 at9.98pt
\begin{document}
\def\semidirect{\hbox{$\;$\bb\char'156$\;$}}
\def\wh#1{\widehat{#1}}
\def\wt#1{\widetilde{#1}}
\newcommand{\SL}{\mathrm{SL}}
\newcommand{\ls}{\mathrm{sl}}
\newcommand{\GL}{\mathrm{GL}}
\newcommand{\g}{{{\mathfrak g}\,}}
\newcommand{\bor}{{{\mathfrak b}}}
\newcommand{\n}{{{\mathfrak n}}}
\newcommand{\h}{{{\mathfrak h\,}}}
\newcommand{\Id}{{\operatorname{Id}}}
\newcommand{\Z}{{\mathbb Z}}
\newcommand{\ZZ}{{\mathbb Z_{>0}}}
\newcommand{\N}{{\mathbb N}}
\newcommand{\R}{{\mathbb R}}
\newcommand{\p}{{\mathbb P}}
\newcommand{\C}{{\mathbb C}}
\newcommand{\Q}{{\mathbb Q}}
\newcommand{\D}{\mathcal{D}}
\newcommand{\I}{\mathcal{I}}
\newcommand{\bF}{{\bf F}}
\newcommand{\LL}{\mathcal{L}}
\newcommand{\F}{\mathcal{F}}
\newcommand{\W}{\mathcal{W}}
\newcommand{\PP}{\mathcal{P}}
\newcommand{\KK}{\mathcal{K}}
\newcommand{\Sym}{{\rm Sym}}
\newcommand{\Sing}{{\rm Sing}}
\newcommand{\Poly}{{\rm Poly}}
\newcommand{\Span}{{\rm Span}}
\newcommand{\Res}{{\rm Res}}
\newcommand{\1}{{\bf 1}}
\newcommand{\A}{{\bf a}}
\newcommand{\m}{{\bf m}}
\newcommand{\z}{{\bf z}}
\newcommand{\uu}{{\bf u}}
\newcommand{\vv}{{\bf v}}
\newcommand{\w}{{\bf w}}
\newcommand{\kk}{{\bf k}}
\newcommand{\T}{{\bf t}}
\newcommand{\s}{{\bf s}}
\newcommand{\dontprint}[1]
{\relax}
\newtheorem%
{thm}{Theorem}
\newtheorem%
{prop}
{Proposition}
\newtheorem%
{lemma}
{Lemma}
\newtheorem%
{lemmadef}[thm]{Lemma-Definition}
\newtheorem%
{cor}
{Corollary}
\newtheorem%
{conj}
{Conjecture}

\newtheorem%
{remark}
{Remark}

\newtheorem%
{definition}
{Definition}

\title {Intersections of Schubert varieties and\\
highest weight vectors \\
in tensor products of  $sl_{N+1}$-representations}

\author[{}]
{I.~Scherbak}

\begin{abstract}${}$
There is a correspondence between highest weight vectors in the
tensor product of finite-dimensional irreducible
$sl_{N+1}$-modules marked by distinct complex numbers, on the one
hand, and elements of the intersection of the Schubert varieties
taken with respect to the osculating flags of the normal rational
curve at the points corresponding to these complex numbers, on the
other hand, \cite{MV1}, \cite{Sc1}, \cite{Sc3}. The highest weight
vectors are the Bethe vectors of the $sl_{N+1}$ Gaudin model  and
the elements are the $(N+1)$-dimensional non-degenerate planes in
the vector space of complex polynomials.

In the present paper we exploit this correspondence in order to 
calculate Bethe
vectors is the tensor product of two irreducible finite-dimensional
$sl_{N+1}$-representations. 
We find  the Bethe vector  in the case when one of  the two
representations is a symmetric power of the standard one.
The idea is to look for the  intersection of Schubert varieties 
related to a Bethe vector.
We present explicitly  a basis of the corresponding
$(N+1)$-dimensional plane in the space of polynomials.

\end{abstract}
\maketitle
\pagestyle{myheadings}
\markboth{I.~Scherbak}
{Intersections of Schubert varieties and highest weight vectors}

\medskip

\section{Introduction}\label{s1}

A {\it partition} $\w=(w_1,\dots,w_N)$ is a collection of non-negative
integers in weakly decreasing order,
$w_1\geq w_2\geq \dots \geq w_N$.
A partition with at most one non-zero entry is called {\it special} and
denoted by $(w_1)$.

\medskip
As it is well-known (\cite{Fu}), any partition $\w$
defines the finite-dimensional irreducible $sl_{N+1}$-representation
$L_\w$ with highest weight
$\Lambda_\w=w_1\lambda_1+\dots+w_N\lambda_N\in\h^*$, where
$$
\h^*\,=\,\C\{\lambda_1,\dots,\lambda_{N+1}\}\,/\,
\left( \lambda_1+\dots+\lambda_{N+1}=0 \right)\,
$$
is the dual to the Cartan subalgebra. The same partition
determines also a Schubert variety $\Omega_\w$ of the
complex codimension
$$
|\w|=w_1+\dots\,+w_N
$$
in the Grassmannian of $(N+1)$-dimensional linear subspaces of $\C^{d+1}$,
where $d\geq N+w_1$ and $w_{N+1}:=0$.
If a partition is special, $\w=(m)$, then $L_{(m)}$ is
the $m$-th symmetric power of the standard $sl_{N+1}$-representation,
and $\Omega_{(m)}$ is a special Schubert variety.

\medskip
This correspondence between highest weight representations and
Schubert varieties can be extended to a correspondence between
certain highest weight vectors in the tensor product of
finite-dimensional irreducible $sl_{N+1}$-modules marked by
distinct complex numbers, on the one hand, and certain elements in
the intersection of the Schubert varieties taken with respect to
the osculating flags of the normal rational curve at the points
corresponding to these complex numbers, on the other hand,
\cite[Section 5]{MV1}, \cite{Sc1} (the case $N=1$ has been done in
\cite{Sc3}).

\medskip
Such highest weight vectors appear in the Gaudin model of
statistical mechanics, as common eigenvectors of certain
mutually commuting linear operators. They are obtained by the Bethe Ansatz
method and called {\it Bethe vectors}, see \cite{FaT},
\cite{FeFR}, \cite{RV} and references therein. The corresponding
elements of the intersection of Schubert varieties are  {\it
non-degenerate planes} in the vector space of complex polynomials,
\cite{Sc1}, \cite{Sc3}.

\medskip
Consider the Grassmannian  $Gr_{N+1}({\Poly}_d)$ of
$(N+1)$-dimensional subspaces of the space of complex polynomials
in one variable of degree at most $d$. We will assume $d$ to be
big enough. For every $\xi\in\C\cup\infty$, denote by
$\F_\bullet(\xi)$ the osculating flag at $\xi$ (that is, the flag
defined by the order of polynomials at $\xi$).

\medskip
Fix partitions $\w(1),\,\dots\,,\,\w(n),\,\w(n+1)$ such that the
sum of codimensions equals the dimension of $Gr_{N+1}({\Poly}_d)$,
$$
|\w(1)|+\dots\,+|\w(n+1)|=(N+1)(N-d)\,.
$$
Let $z_1,\dots,\,z_n$ be distinct complex numbers,
$z=(z_1,\dots,\,z_n)$. Consider the intersection of Schubert
varieties taken with respect to the osculating flags,
\begin{equation}\label{In}
\I_{\{\w\}}(z)\,:=\,\Omega_{\w(1)}\left(\F_\bullet(z_1)\right)\cap
\dots \cap \Omega_{\w(n)}\left(\F_\bullet(z_n)\right)\cap
\Omega_{\w(n+1)}\left(\F_\bullet(\infty)\right) \subset
G_{N+1}({\Poly}_d)\,,
\end{equation}
and the tensor product of $n$
irreducible finite-dimensional representations,
\begin{equation}\label{L}
L\,:=\,L_{\Lambda_{\w(1)}}\otimes\,\dots\,\otimes L_{\Lambda_{\w(n)}}\,,
\end{equation}
marked by $z_1,\dots,\,z_n$, respectively. Denote by $\w^*(n+1)$
the partition dual to $\w(n+1)$.

\begin{thm} {\rm  (\cite{MV1}, \cite{Sc1} -- \cite{Sc4})}
There is a one-to-one correspondence between the non-degenerate
planes in   $\I_{\{\w\}}(z)$\, and the Bethe vectors of the weight
$\Lambda_{\w^*(n+1)}$ in the Gaudin model associated with $L$ and
$z$.
\end{thm}
 Definitions of  Bethe vectors and
non-degenerate planes are given in Sections~\ref{s24}
and~\ref{s34}, respectively.

\medskip
The correspondence, which is crucial in the present work, is given by
a remarkable symmetric
rational function  called {\it the master function of the model}
in the Bethe Ansatz (\cite{ScV}, \cite{MV1}) and {\it the
generating function of the Schubert intersection} in the Schubert
calculus (\cite{Sc2}--\cite{Sc4}).
Both the Bethe vectors and the
non-degenerate planes are determined by the orbits
of critical points with non-zero critical value of this function.

\begin{conj}\label{Cj1} {\rm (The Schubert calculus conjecture)}\ \
For generic $z=(z_1,\dots\,,z_n)$ the intersection of Schubert varieties
 $\I_{\{\w\}}(z)$ is transversal.
\end{conj}

The words ``generic $z$'' mean ``$z$ does not belong to a suitable
proper algebraic surface in $\C^n$''. In the case when at least
$(n-1)$ of the partitions $\w(1),\dots\,,\w(n+1)$ are special, the
transversality has been proved in~\cite{EH}, \cite{So}.

\begin{conj}\label{Cj3}{\rm (Non-triviality of Bethe vectors, \cite{CSc})}\ \
 In the $sl_{N+1}(\C)$ Gaudin model,  every Bethe vector is non-zero,
 for any $z$.
\end{conj}

This conjecture holds for $N=1$  (\cite{ScV}, \cite{Sc2}).
For $N>1$ it was verified in certain
examples in~\cite{MV2} and \cite{CSc}.

\medskip
Well-known relations between the Schubert calculus and
representation theory via the Littlewood--Richardson coefficients
(\cite{Fu}; see also Proposition~\ref{PFu} in our
Section~\ref{s36}) and recently established relations between
Gaudin model and opers (\cite[Section~5]{F}) lead to the following
conclusions.

\medskip\noindent
{\bf Claim.}\ \

\begin{itemize}
\item {\it Assume that  {\rm Conjecture~\ref{Cj3}} holds. Assume that for
some $z^{(0)}$ the intersection {\rm (\ref{In})} with $z=z^{(0)}$
is transversal and consists of non-degenerate planes only.  Then
for generic $z$, the Bethe vectors of the weight
$\Lambda_{\w^*(n+1)}$ in {\rm (\ref{L})} form a basis in the
subspace of singular vectors of this weight.}
\item {\it Assume that for some $z^{(0)}$ the Bethe vectors form
a basis in the subspace of singular vectors  of the weight
$\Lambda_{\w^*(n+1)}$ in the Gaudin model associated with {\rm
(\ref{L})} and $z=z^{(0)}$. Then for generic $z$ the intersection
{\rm (\ref{In})} is transversal and consists of non-degenerate
planes only.}
\end{itemize}

\bigskip
In the Gaudin model, the case of the tensor product of $n=2$
irreducible finite-dimensional representations is  basic. In
the basic case all values of $z$ are generic (so one can assume,
say, $z=(0,1)$) and the operators are reduced to the Casimir
operator. The existence and properties of Bethe vectors for $n>2$
could be deduced from the existence of Bethe vectors in certain
basic cases, applying the construction of  iterated singular
vectors presented in~\cite[Sec.~8, 10]{RV} and results
of~\cite[Section~5]{F}.

\begin{thm}
If $L$ and $\w^*(n+1)$
satisfy the {\rm Basic Property}  , then for generic $z$
there exists a basis of the Bethe vectors in the corresponding Gaudin model.
\end{thm}
The Basic Property is formulated in Section~\ref{s45}.
We refer to the existence of a basis of Bethe vectors as to {\it the Bethe Ansatz
completeness}.

\medskip
In the Schubert calculus in $Gr_{N+1}({\Poly}_d)$,  the
intersection of  $3(=2+1)$ Schubert varieties is basic. Again, for
the basic intersection all values of $z$ are generic.
According to B.~Osserman, transversality  of
every basic intersection implies  transversality
of the intersection of $n+1>3$ Schubert varieties,  for generic $z$
(\cite[Theorem~1.3]{O}).

\medskip
In the proofs of the Bethe Ansatz completeness for $N=1$  (\cite{ScV})) and for
the tensor product of several copies of first and last fundamental
$sl_{N+1}$-modules  (\cite{MV2}), to find  Bethe vectors in the basic case
was the initial and rather sophisticated  step.
Relations to the Schubert calculus suggest to look for
non-degenerate planes in place of Bethe vectors. 

\medskip
In this paper
first steps in this direction are done. We calculate explicitly
the (single element of the) intersection of three Schubert
varieties in the case when at least one of the partitions is
special.

\medskip
In particular, we obtain the following result. Call the polynomial
$$
P_{m;\,d}(x):=1+dx+{d\choose 2}x^2+\dots +{d\choose m}x^m\,,\ \
m\leq d\,,
$$
a {\it truncated binomial}. In particular, $P_{d;\,d}(x)=(x+1)^d$.

\begin{thm} The truncated binomials
$$
P_{m_1;\,d}(x)\,,\ \dots\,,\ P_{m_N;\,d}(x)\,,\ P_{d;\,d}(x)\,,
$$
where $0\leq m_1<m_2<\dots<m_N<d$, span (the single element of)
the basic intersection
$$
\I_{\A,\,\w}\,=\,\Omega_{\A}\left(\F_\bullet(0)\right)\cap
\Omega_{(d-N)}\left(\F_\bullet(-1)\right)\cap
\Omega_{\w}\left(\F_\bullet(\infty)\right)\,\subset\,Gr_{N+1}(\Poly_d)\,,
$$
where
$$
\begin{array}{l}
\A\,=\,(m_N+1-N, \dots\,, m_2-1,m_1 )\,,\\
\w \,=\,(d-m_1-N, \dots\,,d-m_{N-1}-2, d-m_N-1)\,.
\end{array}
$$
\end{thm}

The generic case is done in Theorem~\ref{T5}. If the obtained
element is a non-degenerate plane, then our calculation gives also
a Bethe vector in the tensor product of two irreducible
$sl_{N+1}$-representations where at least one of them is a
symmetric power of the standard one. For $N=1$ the Bethe vector
was calculated explicitly in~\cite{V}, by another method.

\medskip\noindent{\bf Plan of the paper}\ \
In Section~\ref{s2} we collect data on Bethe vectors in the
$sl_{N+1}$ Gaudin model. In Section~\ref{s3} we describe  Schubert
intersections and the generating function, following~\cite{Sc1},
\cite{Sc4}. In Section~\ref{s36} relations between non-degenerate
planes and Bethe vectors are explained. Section~\ref{s4} is
devoted to the basic case, and in Section~\ref{s5} we calculate
the intersection of three Schubert varieties when one of them is
special. If the obtained intersection is non-degenerate, we get
the corresponding Bethe vector.

\bigskip\noindent
{\bf Acknowledgments.}\ \
This work has been done in February--June 2004, when the author visited
the MSRI and the Mathematical Department of the Ohio State University.
Results of  Sections~\ref{s43} and  \ref{s5} have been presented
 in the seminar of D.~Eisenbud on algebraic geometry
at the University of California, Berkeley, on February 2004.
I am  grateful to these institutions for hospitality and excellent working
conditions. It is my pleasure to thank S.~Chmutov and E.~Frenkel for many
useful discussions.

\section{Bethe vectors}\label{s2}
\subsection{Finite-dimensional irreducible $sl_{N+1}$-modules}\label{s21}
For the basic notions of representation theory see for example~\cite{FuH}.

\medskip
Denote by $\{e_i,f_i,h_i\}_{i=1}^{N}$
the standard Chevalley generators of the Lie algebra
$sl_{N+1}(\C)$,
$$
[h_i,e_i]=2e_i\,,\ [h_i,f_i]=-2f_i\,,\ [e_i,f_i]=h_i\,;\qquad
[h_i,h_j]=0\,,\ [e_i,f_j]=0\,  {\rm \ if \ } \ i\neq j\,;
$$
denote by $\h^*$ the dual to the Cartan subalgebra,
$$
\h^*\,=\,\C\{\lambda_1,\dots,\lambda_{N+1}\}\,/\,
\left( \lambda_1+\dots+\lambda_{N+1}=0 \right)\,,
$$
equipped with the standard bilinear form $(\cdot,\cdot)$.
The simple positive roots are  $\alpha_i=\lambda_i-\lambda_{i+1}$,
$1\leq i\leq N$,
$$
(\alpha_i,\alpha_i)=2\,;\quad (\alpha_i,\alpha_{i\pm 1})=-1\,;
\quad (\alpha_i,\alpha_j)=0\,, \mbox{\ if\ } |i-j|>1\,.
$$

\medskip
A partition $\w=(w_1,\dots,w_N)$
defines the finite-dimensional irreducible $sl_{N+1}$-representation
$L_\w$ with highest weight
$\Lambda_\w=w_1\lambda_1+\dots+w_N\lambda_N\in\h^*$.

\subsection{The Gaudin model \cite{G}, \cite{F}}\label{s22}
Fix $n$ partitions $\w(1),\dots\,,\w(n)$ and consider the
tensor product
$$L=L_{\w(1)}\otimes\dots\otimes  L_{\w(n)}\, ,
$$
where $L_{\w(j)}$  is a finite-dimensional irreducible
$sl_{N+1}$-module  with the highest weight
$$\Lambda_{\w(j)}\,=\,w_1(j)\lambda_1+\dots+
w_N(j)\lambda_N\,,\quad 1\leq j\leq n\,.
$$
In  the Gaudin model of statistical mechanics, $L_{\w(1)}, \dots,
L_{\w(n)}$  are labelled by distinct complex numbers
$z_1,\dots,z_n$, and $L$ is called the space of states of the
model. Write $z=(z_1,\dots,z_n)$. Bethe vectors are common
eigenvectors of the mutually commuting linear operators
$H_1(z),\dots\,, H_n(z)$ in $L$ which are defined as follows,
\begin{equation}\label{H}
H_j(z)=\sum_{i\neq j}\, \frac{C_{ij}}{z_j-z_i}\,,\ \
1\leq j\leq n\ ,
\end{equation}
here $C_{ij}$ acts as the Casimir operator on factors
$L_{\w(i)}$ and
$L_{\w(j)}$ of the tensor product and as the identity
on all other factors.

\medskip
The main problem in the Gaudin model is to find a common eigenbasis
and the spectrum of $H_1(z),\dots\,, H_n(z)$.
The operators commute with the diagonal
action of $sl_{N+1}$ in $L$, therefore it is enough to construct common
eigenvectors in the subspace of singular vectors of a given weight,
for every weight.

\medskip
Choose highest weight vectors
$\vv_j\in L_{\w(j)}$,\ $j=1,\dots,n$. Clearly
$\vv_1\otimes\dots\otimes \vv_n\in L$
has the maximal weight and is a common eigenvector
of  $H_j(z)$'s. The idea of the {\it Bethe Ansatz} is to construct
eigenvectors of other weights
by applying certain operators (depending on auxiliary parameters)
to vectors $\vv_1\,,\dots\,,\, \vv_n$.
This idea was realized
in~\cite[Sections~6,7]{SV} for any simple Lie algebra, in the context
of the Knizhnik--Zamolodchikov equation, see also~\cite{BaFl}.
We explain the construction of V.~Schechtman and A.~Varchenko
in the next subsection.

\subsection{The universal weight function} \label{s23}

Highest weight vectors in the tensor product $L$ have weights of the form
\begin{equation}\label{Lambda}
\Lambda(k)\,=\,\Lambda(1)+\dots +\Lambda(n)
-k_1\alpha_1-\dots-k_N\alpha_N\,,
\end{equation}
where $k_1,\dots,k_N$
are non-negative integers  such that $(\Lambda(k)\,,\,\alpha_j)\geq 0$
for every $1\leq j\leq N$. Fix such $k=(k_1,\dots,k_N)$.
We shall construct a function $\vv(\T)$ depending on some
auxiliary variables $\T$ and taking values in the weight subspace
$L_{\Lambda(k)}\subset L$ of weight $\Lambda(k)$.
As we will see in the next subsection, for certain values of the auxiliary
variables, the values of this {\it universal weight function}
will be common eigenvectors of the Gaudin operators $H_j(z)$'s.
The universal weight function is constructed in four steps described
below.

\medskip\noindent
\underline{\it Step I}\ \ is to
choose vectors that  generate the weight subspace
$L_{\Lambda(k)}$.

Consider all  $n$-tuples of words $(\bF_1,\dots, \bF_n)$ in
letters $f_1,\dots, f_N$ subject to the condition that the total
number of occurrences of letter $f_i$ is precisely $k_i$. Our
vectors will be labelled by these $n$-tuples.

Namely, we may think about  $\bF_i$ as an element of the universal
enveloping algebra of $sl_{N+1}$ that naturally acts on the space
$L_{\w(i)}$,\ $i=1,\dots,n$. Then the vector

\begin{equation}\label{u}
 \uu_{_{(\bF_1,\dots\,,\,\bF_n)}} := \bF_1\vv_1\,\otimes\,\dots\,
\otimes\, \bF_n\vv_n
\end{equation}
has weight $\Lambda(k)$, and all such vectors generate the weight
subspace $L_{\Lambda(k)}\,$. In general, their number is greater than the
dimension of that subspace, so they are linearly dependent.

\medskip\noindent
\underline{\it Step II}\ \ For every $i=1,\dots, N$,
introduce a set of $k_i$ auxiliary variables associated with
the root $\alpha_i$,
$$
t(i):= \left(\, t_1(i)\,,\  \dots\,,\  t_{k_i}(i)\right)\,,
$$
and  write $\T:=\left(\, t(1)\,,\ \dots\,, t(N)\,\right)$.

We define  $\vv_{k,\,z}(\T)$ as a linear combination
of the vectors constructed in the first step,

\begin{equation}\label{v}
\vv_{k,\,z}(\T)\ :=\ \sum_{(\bF_1,\dots\,,\,\bF_n)}\
      \omega_{_{(\bF_1,\,z_1,\dots\,,\,\bF_n, \,z_n)}}(\T)\,
\uu_{_{(\bF_1,\dots\,,\,\bF_n)}}\,,
\end{equation}
where $\omega_{_{(\bF_1,z_1,\dots\,,\, \bF_n,z_n)}}(\T)$ are certain
rational functions. These functions are constructed
in the two next steps.

\medskip\noindent
\underline{\it Step III} \ \
Let us write-down the words
$$\bF_1 = f_{i_{_{1,1}}}\dots f_{i_{_{1,s_1}}}\,;\ \dots\,;\
  \bF_n = f_{i_{_{n,1}}}\dots f_{i_{_{n,s_n}}}\,;\ \
  \bF_1\,\dots\,\bF_n = f_{i_{_{1,1}}}\dots f_{i_{_{1,s_1}}}\,
\dots\, f_{i_{_{n,1}}}\dots f_{i_{_{n,s_n}}}\,.
$$
The length of the word $\bF_1\dots\bF_n$ equals
$s_1+\dots +s_n=k_1+\dots+k_N$.

Now we translate $\bF_1\,\dots\, \bF_n$ and
$z_1,\dots\,,\,z_n$
into a rational function
$g_{_{(\bF_1,z_1,\dots\,,\,\bF_n,z_n)}}(\T)$ of $\T$ in the
following way.
For every  $i=1,\dots, N$, we replace the
first occurrence (from left to right) of $f_i$ in the world
$\bF_1\dots \bF_n$ by the variable $t_1(i)$; the second occurrence
by the variable $t_2(i)$; and so on up to the last, $k_i$-th,
occurrence, where $f_i$ will be replaced by $t_{k_i}(i)$. We will get a
certain $n$-tuple of words in $\T$.  Augmenting the $j$-th word in this
$n$-tuple by $z_j$, we get the row,
$$t_{a_{_{1,1}}}(i_{_{1,1}})\,,\, t_{a_{_{1,2}}}(i_{_{1,2}})\,,\,\dots\,
\,,\, t_{a_{_{1,s_1}}}(i_{_{1,s_1}})\,,\, z_1\,,
  \dots\,,
  t_{a_{_{n,1}}}(i_{_{n,1}})\,,\, t_{a_{_{n,2}}}(i_{_{n,2}})\,,\,\dots\,,\,
t_{a_{_{n,s_n}}}(i_{_{n,s_n}})\,,\, z_n\,,
$$
in which every variable $t_a(i)$ from $\T$ appears precisely once.
This row defines the product of fractions

\begin{equation}\label{g}
\begin{array}{l} \displaystyle
g_{_{(\bF_1,z_1,\dots,\,\bF_n,z_n)}}(\T) \,:=\,

\prod_{j=1}^{s_1-1}\frac1{\left(t_{a_{_{1,j}}}(i_{_{1,j}})-
t_{a_{_{1,j+1}}}(i_{_{1,j+1}})\right)}\times\,
\frac 1{\left(t_{a_{_{1,s_1}}}(i_{_{1,s_1}})-z_1\right)}\
\ldots \vspace{15pt} \\ \displaystyle \ldots\
\prod_{j=1}^{s_n-1}\frac1{\left(t_{a_{_{n,j}}}(i_{_{n,j}})-
t_{a_{_{n,j+1}}}(i_{_{n,j+1}})\right)}\times\,
\frac 1{\left(t_{a_{_{n,s_n}}}(i_{_{n,s_n}})-z_n\right)}\,.
\end{array}
\end{equation}

\medskip\noindent
\underline{\it Step IV}\ \ is  symmetrization of
 $g_{_{(\bF_1,z_1,\dots,\,\bF_n, z_n)}}(\T)$.
Let $S_{(k)}$ denote the group of permutations of variables
$$\T=\left(t_1(1),\dots, t_{k_1}(1),\ t_1(2),\dots,\ t_{k_2}(2),
    \ \dots\ ,\ t_1(N),\dots,\ t_{k_N}(N)\right)
$$
that permute variables $t_1(i),\dots,\,t_{k_i}(i)$ within their own,
$i$-th, set, for every $i=1,\dots,N.$
Thus $S_{(k)}$ is isomorphic to the direct product
$S_{k_1}\times S_{k_2}\times\dots\times S_{k_N}$ of $N$ complete
groups of permutations.

For a function $g(\T)$ define the symmetrization by the formula
$$\s_{(k)}[g] :=  \sum_{\pi\in S_{(k)}}
g\Bigl(\pi\bigl(t_1(1)\bigr),\dots, \pi\bigl(t_{k_N}(N)\bigr)\Bigr)\,.
$$

Finally we set
\begin{equation}\label{omega}
\omega_{_{(\bF_1,\,z_1,\dots\,,\, \bF_n,\,z_n)}}(\T) :=
\s_{(k)}\left[g_{_{(\bF_1,\,z_1,\dots\,,\bF_n,\,z_n)}}(\T)\right]\,.
\end{equation}
The function $\vv(\T)=\vv_{k,z}(\T)$  given by (\ref{u}) --
(\ref{omega}) was called in~\cite{MV2}
{\it the universal weight function}. It takes values in $L_{\Lambda(k)}$.

\subsection{Bethe vectors of the weight $\Lambda(k)$}\label{s24}

For convenience, in what follows we set $w_{N+1}(j)=0$ for $1\leq j\leq n$
and $k_0=k_{N+1}=0$.
\begin{thm}\label{T1} {\rm \cite{BaFl}, \cite{FeFR}, \cite{RV}}
Vector $\vv(\T)$ given by {\rm (\ref{u}) -- (\ref{omega})}
is a common eigenvector of
$H_1(z),\dots, H_n(z)$ if and only if {\rm the Bethe equations}
$$
-\,\sum_{j=1}^n\ \frac{\left(\alpha_i,\Lambda_{\w(j)}\right)}
{t_l(i)-z_j}\,+\,\sum_{s\neq l}\ \frac 2{t_l(i)-t_s(i)}\,-\,
\sum_{s=1}^{k_{i-1}}\ \frac 1{t_l(i)-t_s(i-1)}\,-\,
\sum_{s=1}^{k_{i+1}}\ \frac 1{t_l(i)-t_s(i+1)}\,=0\,,
$$
where $1\leq i\leq N$, $1\leq l\leq k_i$,  are satisfied.
Moreover, if $\vv(\T)$ is an eigenvector,
then it is a highest weight vector.
\end{thm}

\begin{remark} \label{r1}{\rm For $1\leq j\leq n$, we have
$$
\left(\alpha_i,\Lambda_{\w(j)}\right)=w_i(j)-w_{i+1}(j)\,,\ \
1\leq i\leq N\,.
$$
The weight $\Lambda(k)$ given by  (\ref{Lambda}) corresponds to the partition
$\w(k)=\left(w_1(k),,\, \dots\,,\, w_N(k)\right)$, where
$$
w_i(k)=\sum_{j=1}^n\,w_i(j)\,-\,k_i\,+\,k_{i-1}-k_N\,,
\quad 1\leq i\leq N\,,
$$
and thus we write $\Lambda(k)\,=\,\Lambda_{\w(k)}$.}
\end{remark}

\begin{definition} \label{d1}
{\rm The value of the universal weight function $\vv(\T)$ at
a solution of the Bethe equations is called a}  Bethe vector.
\end{definition}

\subsection{Master function}\label{s25}
The function
$$
\begin{array}{l}\vspace{15pt}\displaystyle
\Psi(\T)\,=\,\prod_{j=1}^n\prod_{i=1}^N\prod_{l=1}^{k_i}
 \left(t_l(i)-z_j\right)^{-(\alpha_i,\,\Lambda_{\w(j)})}
\prod_{i=1}^N\prod_{1\leq l<s\leq k_i}
 \left(t_l(i)-t_s(i)\right)^2\,\times\\ \vspace{15pt}\displaystyle
\quad\quad \quad \times\,\prod_{i=1}^{N-1}\prod_{l=1}^{k_i}
\prod_{s=1}^{k_{i+1}}\left(t_l(i)-t_s(i+1)\right)^{-1}
\end{array}
$$
is called the {\it master function} associated with
$L_{\Lambda(k)}\subset L$ and $z$.

As it was pointed out in~\cite{RV},

\begin{itemize}
\item  the Bethe equations of Theorem~\ref{T1} are exactly
the defining equations of the critical points with non-zero critical
values of the function $\Psi(\T)$;
\item
this function is symmetric with respect to the action of group $S_{(k)}$
defined in Step IV, Section~\ref{s23},
and the critical points belonging to the same orbit  define
the same Bethe vector.
\end{itemize}

On critical points of the master function see also
 \cite{MV1},  \cite{MV2}, \cite{ScV}, \cite{Sc1}.

\section{Non-degenerate planes in Schubert intersections}\label{s3}

\subsection{Schubert cells in the Grassmannian of
$(N+1)$-dimensional planes in $\Poly_d$}\label{s31}
For the Schubert calculus see, for example,~\cite{Fu}.

\medskip
We identify $\C\p^d$ endowed with the embedded normal rational curve,
$$
x\mapsto [1:x:x^2:\dots :x^d]\,,\quad
\infty\mapsto [0:0:\dots 0:1]\,,\quad  x\in \C\,,
$$
and the vector space $\Poly_d$ of complex polynomials in
$x$ of degree at most $d$, considered up to a non-zero factor.

\medskip
Fix $\xi\in\C\cup\infty$ and take the  flag $\F_\bullet(\xi)=
\{\F_0(\xi) \subset \F_1(\xi)\subset \dots\subset \F_{d}(\xi)=
{\Poly}_{d}\}$ defined by the order of polynomials at $\xi$. That
is, if $\xi\in\C$, then $\F_i(\xi)$ consists of the polynomials of
the form $a_i(x-\xi)^{d-i}+\dots +a_{0}(x-\xi)^{d}$, and
$\F_i(\infty)={\Poly}_i\,$, $0\leq i\leq d$ (by definition,
the order of $f(x)\in{\Poly}_d$ at $\infty$ is $d-\deg f(x)$).
Thus
$\F_\bullet(\xi)$ is the osculating flag of the normal rational
curve in $\C\p^d$ at $\xi$.

\medskip
Denote by $Gr_{N+1}({\Poly}_d)$ the Grassmannian of $(N+1)$-dimensional
subspaces of $\Poly_d$. This is an algebraic variety of dimension
$(d-N)(N+1)$.

\medskip
The integers $w_1,\dots, w_{N+1}$
determine the {\it Schubert cell} with respect to flag  $\F_\bullet(\xi)$.
This cell is formed  by the elements  $V\in  Gr_{N+1}({\Poly}_d)$ satisfying
the conditions
$$
\dim \left(V\cap \F_{d-N-1+i-w_i}(\xi)\right)= i\,,\quad \dim
\left(V\cap \F_{d-N-1+i-w_i-1}(\xi)\right)= i-1\,, \quad
i=1,\dots, N+1\,.
$$
Clearly the Schubert cell  may be non-empty only if
$d-N\geq w_1\geq \dots w_{N+1}\geq 0\,.$ Notice also that if $w_{N+1}>0$,
then all polynomials in $V$ are divisible by $(x-\xi)^{w_{N+1}}$, that is
$V$ {\it has a base point} at $\xi$.

\subsection{Schubert intersections}\label{s32}

In what follows, we will consider Schubert cells such
that  $w_{N+1}=0$ and will denote  by
$$\Omega^\circ_{\w}(\xi)\subset Gr_{N+1}({\Poly}_d)$$
the Schubert cell corresponding to the flag $\F_\bullet(\xi)$ and to
the partition $\w=(w_1,\dots\,,w_N)$ with
$d-N\geq w_1\geq \dots w_N\geq 0\,$. In other words, we  will consider
elements $V\in Gr_{N+1}({\Poly_d})$ {\it with no base point}.
This means that for any $\xi\in\C$ there
is a polynomial in $V$ which does not vanish at $\xi$, and $V$
contains a polynomial of degree $d$ (by definition, $f(x)\in\Poly_d$
has a root of multiplicity $d-\deg f(x)$ at $\infty$; so every polynomial
in $\Poly_d$ has exactly $d$ roots counting with multiplicities).
This assumption is not very restrictive in fact; it says that $d$ is
chosen as small as possible and that the polynomials of $V$ do not have
any common factor.

\medskip
The number  $|\w|=w_1+\dots\,+w_N$ is the complex
codimension of $\Omega^\circ_{\w}(\xi)$ in  $Gr_{N+1}(\Poly_d)$.
The closure of a Schubert cell in the Grassmannian,
$\Omega_\w(\xi)=\overline{\Omega^\circ}_\w(\xi)\,,$
is called  a {\it Schubert variety}. In particular,
$\Omega_{(0,\dots,0)}(\xi)=Gr_{N+1}({\Poly_d})$.
The cohomology class $\sigma_{\w}$
of $\Omega_{\w}(\xi)$ does not depend on flag choice and
is called  the  {\it Schubert class}, see \cite{Fu}.

\medskip
For every $\xi\in\C\cup\infty$ and every $V\in Gr_{N+1}(\Poly_d)$ with
no base point there exists the unique partition
$\w=\w(\xi;V)$ such that $V\in\Omega^\circ_{\w}(\xi)$.

\medskip
{\it The Wronskian} of $V\in G_{N+1}({\Poly_d})$ is defined as
a monic polynomial $W_V(x)$ which is proportional to the Wronski
determinant of some (and hence, any) basis of $V$.
On relations of the Wronskian to the Schubert calculus see
\cite{ErGa}, \cite{KSo}. In particular the following result can
be easily obtained.

\begin{lemma}\label{L1} For every $\xi\in \C$,
the Wronskian $W_V(x)$ has the order $|\w(\xi;V)|$ at $\xi$,
and the order of $W_V(x)$ at $\infty$ is
$\, \left(\dim Gr_{N+1}({\Poly_d})-\deg W_V\right)$, so that
$$
\sum_{\xi\in\C\cup\infty}\ \ |\w(\xi;V)|=\dim G_{N+1}({\Poly_d}).
$$
If $\xi\in\C$, then $\w(\xi;V)=(0,\dots,0)$ if and only if
$W_V(\xi)\neq 0$, and  $\w(\infty;V)=(0,\dots,0)$ if and only if
$\deg W_V(x)=\dim G_{N+1}({\Poly_d})$, i.e. the degree is the
maximal possible. \hfill $\square$
\end{lemma}

\medskip
The roots of $W_V(x)$  (including $\infty$) are called {\it singular}
(or {\it ramification}) points of $V$. Let $z_1,\dots, z_n$ be the
finite singular points of $V$.
The intersection of the Schubert cells
$\Omega^\circ_{\w(\xi;V)}(\xi)$, for all $\xi\in\C\cup\infty$,
coincides with  the intersection
$$
\Omega_{\w(z_1;V)}(z_1)\cap\dots\cap
\Omega_{\w(z_n;V)}(z_n)\cap
\Omega_{\w(\infty;V)}(\infty)\,
$$
(otherwise the Wronskian would have more roots than its degree).
This intersection  is zero-dimensional,
however  it may contain more than one element; the cardinality of this
intersection is bounded from above by the intersection number of the
corresponding Schubert classes,
$$
\sigma_{\w(z_1;V)}\cdot\, ... \,
\cdot\sigma_{\w(z_n;V)}\cdot\sigma_{\w(\infty;V)}\,.
$$
Any other element of this intersection has the same
Wronskian, so it lies in the preimage of a given polynomial under the
{\it Wronski map} sending $V\in Gr_{N+1}(\Poly_d)$ into
$W_V(x)$. On the Wronski map,  see~\cite{ErGa}, \cite{KSo}, \cite{Sc3}.

\medskip
Fix $\{\w\}=\{\w(1), \dots \,,\w(n),\w(n+1)\}$, the set of
$(n+1)$ partitions,  $\w(j)=(w_1(j),\dots w_N(j))$, such that
\begin{equation}\label{wI}
|\w(1)|\,+\,+\,\dots\,+|\w(n+1)|=\dim  G_{N+1}({\Poly}_{d})\,.
\end{equation}
Let $z_1,\dots\,, z_n$ be distinct complex numbers.
Write $z=(z_1,\dots\,, z_n)$.

\begin{definition}\label{d2} {\rm The intersection of
Schubert varieties
\begin{equation}\label{I}
\I_{\{\w\}}(z)=\Omega_{\w(1)}(z_1)\cap \dots
\cap \Omega_{\w(n)}(z_n)\cap
\Omega_{\w(n+1)}(\infty)\subset G_{N+1}({\Poly}_d)
\end{equation}
is called a }  Schubert intersection.
\end{definition}

\medskip
If $V\in \I_{\{\w\}}(z)$, then $V$ has no base point and
its Wronskian is
\begin{equation}\label{Wr}
W_V(x)\,=\, W_{|\w|,\,z}(x)\,=\,
(x-z_1)^{|\w(1)|}\,\dots\, (x-z_n)^{|\w(n)|}.
\end{equation}

\subsection{Intermediate Wronskians}\label{s33}
Let $\xi\in\C\cup\infty$ and $V\in G_{N+1}({\Poly}_d)$.
The collection of  orders at $\xi$ of all polynomials in $V$ consists
of $(N+1)$ distinct numbers, and the smallest order is always $0$
as $V$ does not have a base point. Clearly $\xi$ is not
singular if and only if the orders at $\xi$ are the minimal
possible, i.e. $0,1,\dots,N$.

\medskip
Consider the  Schubert intersection  $\I_{\{\w\}}(z)$ given by
(\ref{wI}), (\ref{I}).
For any  $V\in \I_{\{\w\}}(z)$,  all finite singular points are
$ z_1,\dots, z_n$ and the Wronskian  is as in (\ref{Wr}).
The non-zero orders $\rho_l(z_j)$ at the finite singular points and
$d-d_l$ at infinity are as follows,
\begin{equation}\label{rho}
\rho_l(z_j)=w_l(j)+N+1-l\,,\quad
d-d_l=w_l(n+1)+N+1-l\,,\quad 1\leq l\leq N,
\  1\leq j\leq n\,.
\end{equation}
Here $d_l$'s are degrees of polynomials in $V$.
We have
\begin{equation}\label{E5}
0\leq d_1<d_2<\dots <d_N<d_{N+1}=d,\quad
0=\rho_{N+1}(z_j)<\rho_{N}(z_j)<\dots <\rho_1(z_j)\leq d.
\end{equation}

\medskip
Denote by $V_\bullet$  the  flag obtained by the  intersection of $V$
and  $\F_\bullet(\infty)$ (recall that $\F_l(\infty)=\Poly_l$),
\begin{equation}\label{V}
V_\bullet\,=\,\left\{\,
V_1\subset V_2 \subset \dots \subset V_{N+1}=V\,\right\},\quad
\dim V_l=l.
\end{equation}
The degrees of polynomials in $V_l$ are $d_1,\dots, d_l$.
Denote by $W_l(x)$ the Wronskian of $V_l$, $1\leq l\leq N+1$.
We call it {\it the $l$-th intermediate Wronskian of } $V$.
In particular, the $(N+1)$-th intermediate Wronskian coincides with
$W_V(x)=W_{|\w|,\,z}(x)$  given by (\ref{Wr}).

\begin{remark}\label{rr}{\rm Every $(N+1)$-dimensional plane $V$ in the
vector space of
polynomials is determined by the set of its intermediate Wronskians
$W_1(x),\dots\,, W_{N+1}(x)$.
Indeed, $V$ is the solution space of the ordinary differential equation of
order $N+1$ (here we set $W_0(x)=1$),
$$
\frac{d}{dx}\,\frac{W^2_N(x)}{W_{N-1}(x)W_{N+1}(x)}\, \dots \,
\frac{d}{dx}\,\frac{W^2_2(x)}{W_3(x)W_1(x)}\,\frac{d}{dx}\,
\frac{W^2_1(x)}{W_2(x)W_0(x)}\,
\frac{d}{dx}\,\frac{u(x)}{W_1(x)}=0\,.
$$
As $(N+1)$ linearly independent solutions, i.e. a basis of $V$, one can take
\begin{eqnarray*}
 u_1(x) & = & W_1(x)\,, \\
 u_2(x)&=&W_1(x)\int^x\frac{W_2W_0}{W_1^2}\,,\\
 u_3(x)&=&W_1(x)\int^x\left(\frac{W_2(\xi)W_0(\xi)}{W_1^2(\xi)}
 \int^{\xi}\frac{W_1W_3}{W_2^2}\right)\,,\\
 \dots & \ & \dots \quad \ \ \quad \dots\quad \ \ \quad \dots \quad \ \ \quad \dots  \\
 u_{N+1}(x) & = &
W_1(x)\int^x\left(\frac{W_2(\xi)W_0(\xi)}{W_1^2(\xi)}\int^{\xi}\left(
\frac{W_1(\tau)W_3(\tau)}{W_2^2(\tau)}\int^{\tau}\dots\int^{\eta}\left.
\frac{W_{N-1}W_{N+1}}{W^2_N}\right.\right)\dots\right)\,,
\end{eqnarray*}
see~\cite[Part VII, Section~5, Problem 62]{PSz} or~\cite{Sc4}. }
\end{remark}

\medskip
Define  polynomials
\begin{equation}\label{Z}
Z_i(x)=\prod_{j=1}^n(x-z_j)^{m_j(i)},\quad 1\leq i\leq N+1\,,
\end{equation}
where
$$
m_j(1)=0\,,\ \ m_j(i) =\,\sum_{l=N+2-i}^{N}w_l(j)\,,\ \
2\leq i\leq N+1\,,\ 1\leq j\leq n\,.
$$
In particular,  $Z_1(x)=1$ and $Z_{N+1}(x)$ coincides with (\ref{Wr}).

\begin{lemma}\label{L2}  \cite{Sc1}, \cite{Sc4}
The ratio $W_i(x)/Z_i(x)$ is a polynomial of degree
\begin{equation}\label{k}
i(d-N)-\sum_{l=1}^{i}w_l(n+1)-
\sum_{j=1}^{n}\sum_{l=N+2-i}^Nw_l(j)\,, \  1\leq i\leq N+1.
\end{equation}
\end{lemma}

\subsection{Non-degenerate planes~\cite{Sc1}, \cite{Sc4}}\label{s34}
Lemma~\ref{L2} shows that the Schubert intersection  $\I_{\{\w\}}(z)$
given by (\ref{wI}), (\ref{I}) determines certain roots of the intermediate
Wronskians, namely the roots of polynomials $Z_i(x)$. Denote
$$
T_{N+1-i}(x):=W_i(x)/Z_i(x)\,,\quad 1\leq i\leq N+1\,.
$$
Every $T_j$ is a polynomial of degree $k_j$ given by (\ref{k}), according to Lemma~\ref{L2}. In particular, $T_0(x)=1$, i.e. $k_0=0$. The roots of  $T_j$
we call {\it the additional roots} of the $(N+1-j)$-th  intermediate Wronskian.
If  $\I_{\{\w\}}(z)$  contains more than one element, then  the intermediate
Wronskians of these elements are differ by the additional roots.

\medskip
We intend to distinguish some elements of Schubert intersections.
Denote by $\Delta(f)$ the discriminant of polynomial $f(x)$ and by
$\Res(f,g)$ the resultant of polynomials $f(x),\, g(x)$.

\begin{definition}\label{d3}
{\rm We call $V\in\I_{\{\w\}}(z)$  a} nondegenerate plane
{\rm  if the polynomials} $T_1(x),\dots\,, T_N(x)$
\begin{itemize}
\item
do not vanish at the finite singular points,
$T_i(z_j)\,\neq\, 0\,,\ \ 1\leq j\leq n$;
\item
do not have multiple roots, $\Delta (T_i)\,\neq\, 0$;
\item
any two neighboring of them do not have common roots,
$\Res(T_i,T_{i\pm 1})\,\neq\,0$.
\end{itemize}
\end{definition}

In other words, the additional roots of a nondegenerate plane are as generic
as possible.

\subsection{Generating function}\label{s35}
The generating function of a Schubert intersection is a rational
function such that its critical points determine the nondegenerate
elements in this intersection. It was defined in~\cite{Sc1}, \cite{Sc4}
as follows.

\medskip
For fixed $z=(z_1,\dots,z_n)$,  any monic polynomial
$f(x)$ can be presented
in a unique way as the product of two monic polynomials $T(x)$ and $Z(x)$ which
satisfy
\begin{equation}\label{TZ}
f(x)=T(x)Z(x), \ \   T(z_j)\neq 0,\ \  Z(x)\neq 0
{\rm\ \ for\ any\ } x\neq z_j, \ \ 1\leq j\leq n.
\end{equation}
Define {\it  the relative discriminant of $f(x)$  with respect to $z$} as
$$\Delta_{z}(f)=\frac{\Delta(f)}{\Delta(Z)}=\Delta(T){\rm Res}^2(Z,T),$$
and {\it the relative resultant   of  $f_i(x)=T_i(x)Z_i(x)$, $i=1,2$,
with respect to $z$} as
$${\rm Res}_{z}(f_1, f_2)=\frac{{\rm Res}(f_1,f_2)}{{\rm Res}(Z_1,Z_2)}=
{\rm Res}(T_1, T_2){\rm Res}(T_1, Z_2){\rm Res}(T_2, Z_1).$$

\medskip
If $V\in \I_{\{\w\}}(z)$ is nondegenerate, then the decomposition
$W_i(x)=T_{N+1-i}(x)Z_i(x)$ is exactly the presentation of $W_i(x)$
in the form (\ref{TZ}).

\begin{definition}\label{d4}
{\rm The function
\begin{equation}\label{Phi-w}
\Phi_{ \{\w\},\,z} (T_1,\dots, T_N)=
\frac {\Delta_z(W_1)\cdots \Delta_z(W_N)}
{{\rm Res}_z(W_1,W_2)\cdots
{\rm Res}_z(W_N,W_{N+1})}
\end{equation}
is  called } the generating function of the Schubert intersection
 $\I_{\{\w\}}(z)$.
\end{definition}

A part of the following theorem was obtained originally
by A.~Gabrielov, \cite{Ga}, in his study of the Wronski map.

\begin{thm}\label{T3} {\rm (\cite{Sc1}, \cite{Sc4}) }
There is a one-to-one correspondence between
the critical points with non-zero critical values of the function
$\Phi_{ \{\w\},z} (T_1,\dots, T_N)$
and the nondegenerate planes in  $\I_{\{\w\}}(z)$.  \hfill $\square$
\end{thm}

Namely, every critical point defines  the intermediate Wronskians,
and hence a nondegenerate plane, see Remark~\ref{rr}.
Vice versa, for every nondegenerate plane
one can calculate the intermediate Wronskians, and the corresponding
polynomials $T_i(x)$ give a critical point of the generating function.

\section{Non-degenerate planes and Bethe vectors}\label{s36}
The partition $\w^*$ {\it dual} to $\w=(w_1,\dots\,,w_N)$
is defined by the formula
$$
\w^*=(w_1,w_1-w_N, w_1-w_{N-1},\dots\,, w_1-w_2)\,.
$$
In Remark~\ref{r1}, the partition $\w(k)$ such that $\Lambda(k)=\Lambda_{\w(k)}$
has been written explicitly. A direct
calculation (based on the formula (\ref{k}) for $k_{N+1-i}$ of Lemma~\ref{L2})
shows the following. 

\begin{lemma}\label{L*} We have $\w(k)=\w^*(n+1)$. \hfill $\square$
\end{lemma}

\medskip
The next claim is a manifestation of the well-known relation of the Schubert
calculus to representation theory via the famous Littlewood--Richardson
coefficients (\cite{Fu}).

\begin{prop} \label{PFu}
The dimension of the subspace of singular vectors of the weight
$\Lambda(k)=\Lambda_{\w^*(n+1)}$ in $L=L_{\w(1)}\otimes\dots\,\otimes L_{\w(n)}$
coincides with the intersection number of the Schubert classes
$\sigma_{\w(1)}\cdot {\rm \ ... \ }\sigma_{\w(n)}
\cdot \sigma_{\w(n+1)}$.  \hfill  $\square$
\end{prop}

\medskip
We try to make this relation more precise.
Let us re-write the  function (\ref{Phi-w}) in terms of
unknown roots of polynomials
$T_1,\dots,T_N$. Recall that if  $f(x)=(x-a_1)\dots (x-a_A)$ and
$g(x)=(x-b_1)\dots(x-b_B)$,  then
$$
\Delta(f)=\prod_{1\leq i<j\leq A}(a_i-a_j)^2\,,\quad
{\rm Res}(f,g)=\prod_{i=1}^A\prod_{j=1}^B(a_i-b_j).
$$

Denote by
$$
t^{(i)}= \left(\, t_1^{(i)}\,,\  \dots\,,\  t_{k_i}^{(i)}\right)\,,
\quad i=1,\dots, N,
$$
the roots of $T_i(x)$, that is the {\it additional} (not prescribed by
$\{\w\}$ and $z$) roots of the Wronskian $W_{N+1-i}(x)$.
Write $\T=\left(\, t^{(1)}\,,\ \dots\,, t^{(N)}\,\right)$.
With this notation, the function (\ref{Phi-w}) becomes a rational function
in $k_1+\dots+k_N$ variables,
\begin{eqnarray}\label{Phi-t}
\Phi_{ \{\w\},\,z}(\T)=\\
&=& \prod_{i=1}^N\ \
\prod_{1\leq l<s\leq k_i}\ \ (t^{(i)}_l-t^{(i)}_s)^2 \nonumber \\
&\times &\prod_{i=1}^{N-1}\ \prod_{l=1}^{k_{i}}\ \prod_{s=1}^{k_{i+1}}
(t^{(i)}_l-t_s^{(i+1)})^{-1} \nonumber\\
&\times &\prod_{i=1}^N\prod_{j=1}^n\ \prod_{l=1}^{k_i}\,
(t^{(i)}_l-z_j)^{w_{i+1}(j)-w_i(j)}\,,\nonumber
\end{eqnarray}
which is nothing but the master function associated
with $z=(z_1,\dots,z_n)$ and $L_{\Lambda(k)}\subset\,L$ of
Section~\ref{s25}.
Summarizing, we arrive at the following conclusion.

\begin{cor}\label{C1} There is a one-to-one correspondence between
the nondegenerate planes in $\I_{\{\w\}}(z)$ given by {\rm (\ref{wI}),
(\ref{I})} and
the Bethe vectors of the Gaudin model associated with
$z$ and $L_{\Lambda(k)}\subset\,L_{\w(1)}\otimes\,\dots\, \otimes L_{\w(n)}$,
where $\Lambda(k)=\Lambda_{\w^*(n+1)}$.
More precisely, every critical point $\T^{(0)}$ with non-zero critical value
of the function $\Phi_{ \{\w\},\,z}(\T)$ determines (the additional roots of) the
Wronskians $W_1,\dots\,,W_N$ of a non-degenerate plane $V=V(\T^{(0)})$ (and hence
$V$ itself, according to {\rm Remark~\ref{rr}}) as well as a Bethe vector
$\vv(\T^{(0)})$ where $\vv(\T)$ is the universal weight function defined in
{\rm Subsection~\ref{s23}}.
\hfill $\square$
\end{cor}

\medskip
Thus the auxiliary variables introduced in  Step II
of Section~\ref{s23} appear to be the additional roots of the
intermediate Wronskians, and
in a certain sense, the intermediate Wronskians $W_1,\dots,\, W_N$
correspond to  simple roots $\alpha_N,\dots,\,\alpha_1$, respectively.

\section{Basic case: $n=2$}\label{s4}

\subsection{Bethe vectors for the Casimir operator}\label{s41}
If $n=2$, then all values of
$z$ are generic. Indeed, as it has been pointed out in~\cite[Section~5]{RV}, for any fixed
 $z=(z_1,z_2)$  with $z_1\neq z_2$ the linear change of variables 
$\tilde\T=(\T-z_1)/(z_2-z_1)$ turns the
Bethe system on $\T$  into the Bethe system on $\tilde\T$ with 
$\tilde z=(0,\, 1)$.
Therefore one can assume
$L=L_{\Lambda(0)}\otimes L_{\Lambda(1)}$, where  $L_{\Lambda(0)}$ is
associated with $z_1=0$ and  $L_{\Lambda(1)}$ with $z_2=1$.

\medskip
According to (\ref{H}),
for $z=(0,1)$ the Gaudin hamiltonians $H_1\,,\,H_2$ satisfy
$H_2=-H_1=C$, where $C$ is the Casimir operator. Furthermore,
$C$ acts in any irreducible submodule of the tensor product by multiplication
by a constant and hence definitely has eigenvectors.

\medskip
We call  $n=2$ {\it the basic case} because one can use Bethe
vectors in the tensor products of $n=2$ irreducible
representations  for the study of Bethe vectors when $n>2$. The
way  is explained in the next subsection.

\subsection{Iterated singular vectors and the Bethe Ansatz}\label{s45}
Consider
$$L\,=\,L_{\Lambda(1)}\otimes\dots\otimes L_{\Lambda(n)}\,,
$$
the tensor product of $n$ irreducible finite-dimensional $sl_{N+1}$-modules.

\medskip
For every
$i=1,\dots, n$, fix a highest weight vector
$\vv_i\,\in\,L_{\Lambda(i)}$. Write $\kk=(k_1,\dots\,,\,k_N)$ and
fix a highest weight
$$
\Lambda(\kk)=\sum_{i=1}^n \,\Lambda(i)\,- k_1\alpha_1 - \dots - k_N\alpha_N
$$
presented in $L$.
Denote by $[\kk]$ a set of presentations of numbers $k_1,\dots,\,k_N$
as the sum of $(n-1)$ non-negative integers,
$$
k_j=k_{1,\,j}\,+\,\dots\, +\,k_{n-1,\,j}\,,\quad  j\,=\,1,\dots,\,N\,,
$$
such that for every $1\leq l\leq n-1$ the weight
\begin{equation}\label{kl}
\Lambda[\kk,\,l]\,=\,\sum_{i=1}^{l+1}\,\Lambda( i)\,-\,\sum_{j=1}^N\,
\left(\sum_{i=1}^{l}\, k_{i,\,j}\right)\alpha_j
\end{equation}
is the highest weight in the tensor product
$$
L_{\Lambda[\kk,\,l-1]}\,\otimes\, L_{\Lambda(l+1)}\,,
$$
clearly $\Lambda[\kk,\,n-1]\,=\,\Lambda(\kk)$. 
We also set  $\Lambda[\kk,\,0]\,:=\, \Lambda(1)\,.)$ 

\medskip
Denote by $\KK$ be the set of all such presentations.
We say  that $\Lambda(\kk)$ and $L$ satisfy the Basic Property
if the following holds.

\medskip\noindent{\bf Basic Property.}\ \ {\it For every  $[\kk]\in\KK$ and
every $1\leq l\leq n-1$,  there exists a basis of  Bethe vectors
in the subspace of singular vectors of the weight
$\Lambda[\kk,\,l]$  in the tensor product
$L_{\Lambda[\kk,\,l-1]}\,\otimes L_{\Lambda(l+1)}$ of two
$sl_{N+1}$-representations, where the first one is labelled by
$z_1=0$ and the second one by $z_2=1$.}

\medskip
\begin{thm}
Let   $\Lambda(\kk)$ and $L$ satisfy the Basic Property. Then for generic
$z_1,\dots,\, z_n$  there exists a basis of Bethe vectors in the Gaudin model
associated with  the weight  $\Lambda(\kk)$ in $L$ and $z=(z_1,\dots,\, z_n)$.
\end{thm}

\medskip\noindent{\bf Proof.}\ \  First of all,
the Basic Property allows to construct a basis  of  {\it iterated Bethe vectors}
in the subspace of singular vectors of the weight $\Lambda(\kk)$ in $L$.
The construction is due
to N.~Reshetikhin and A.~Varchenko, \cite[Section~8, 10]{RV}, see
also~\cite[Section~3.2]{MV2}. Let us  explain it briefly.

\medskip
Fix $[\kk]\,\in\,\KK$. First, for
$L_{\Lambda(1)}\otimes L_{\Lambda(2)}$, $(z_1,z_2)=(0,1)$,
and for the highest weight $\Lambda[\kk,\,1]$,  take  a basis of Bethe vectors
$\{\vv_{[1,2]}^j\}$  in the subspace of singular vectors of this weight in
$L_{\Lambda(1)}\otimes L_{\Lambda(2)}$.
This basis can be written explicitly in terms of $\vv_1, \vv_2$ and
solutions to the corresponding Bethe system, according to Theorem~\ref{T1}.

Next, if we decompose $L_{\Lambda(1)}\otimes L_{\Lambda(2)}$ into the direct
sum of irreducible representations,
 then  $L_{\Lambda(1)}\otimes L_{\Lambda(2)}\otimes L_{\Lambda(3)}$
can be written as the direct sum of the tensor products of two irreducible representations,
the first one comes from the decomposition and the second one is $ L_{\Lambda(3)}$.
For every tensor product  of the form
$L_{\Lambda[\kk,1]}\otimes  L_{\Lambda(3)}$
(which is again the tensor product of two irreducible
$sl_{N+1}$-modules with fixed highest weight vectors
$\vv_{[1,\,2]}^j$ and $\vv_3$)  and for  the highest weight
$$
\Lambda[\kk,\,2]=\Lambda[\kk,\,1]+\Lambda(3)-k_{2,\,1}\alpha_1-\dots
-k_{2,\,N}\alpha_N
$$
in this tensor product, we take the basis of Bethe vectors for the
subspace of singular vectors of this weight and $(z_1,z_2)=(0,1)$;
here $L_{\Lambda[\kk,1]}$ is marked by $z_1=0$ and
$L_{\Lambda(3)}$ is marked by $z_2=1$. This basis can be written
explicitly in terms of $\vv_{[1,\,2]}^j$, $\vv_3$ and the
solutions to the relevant Bethe system (that is, in terms of
$\vv_1, \vv_2, \vv_3$ and the solutions of  two Bethe systems
corresponding to certain basic cases). If we continue this
procedure (i.e., decompose $L_{\Lambda[\kk,1]}\otimes
L_{\Lambda(3)}$ into the direct sum of irreducible
representations, write $L_{\Lambda[\kk,1]}\otimes
L_{\Lambda(3)}\otimes  L_{\Lambda(4)}$ as the direct sum of the
tensor products of two irreducible representations, one coming
from the decomposition and the other being $L_{\Lambda(4)}$,
etc.),  then as a result we get certain singular vectors of the
weight $\Lambda(\kk)$ in $L$, written in terms of the vectors
$\vv_i$'s  and solutions of $(n-1)$ Bethe systems corresponding to
basic cases.  These vectors are called {\it iterated Bethe
vectors}.

\medskip
Performing this procedure for every $[\kk]\,\in\,\KK$, we obtain  the basis
of iterated Bethe vectors in the subspace of  singular vectors of the weight
$\Lambda(\kk)$ in $L$.

\medskip
The next claim  is similar to~\cite[Theorem~9.16]{RV} and
\cite[Lemma~4.4]{MV2}. For the reader convenience, we sketch a
proof.

\begin{lemma}\label{itB}
Let   $\Lambda(\kk)$ and $L$ satisfy the Basic Property and
$s$ be a big real number.
Then for every iterated Bethe vector $\vv_B$ of the weight $\Lambda(\kk)$
 in $L$, there exist an integer $\nu$ and a  Bethe
vector $\vv^{(s)}$  of the Gaudin model associated with
$L_{\Lambda(\kk)}\,\subset\, L$ and $z=(s,s^2,\dots, s^n)$,  such that
$$
\vv^{(s)}=s^\nu\left(\vv_B+O(s^{-1})\right)\,, \quad s\to+\infty\,.
$$
\end{lemma}

\bigskip\noindent{\bf Proof of the Lemma.} \ \ Every $[\kk]\,\in\,\KK$ determines
$(n-1)$ weights $\Lambda[\kk,\,l]$, $l=1,\dots, n-1$, corresponding to the
$(n-1)$ iteration steps, see (\ref{kl}).

 Let $\Psi=\Psi(\T)$ be the master function of the Gaudin model associated
with $L_{\Lambda(\kk)}\,\subset\, L$ and $z=(s,s^2,\dots, s^n)$.
In what follows we do not distinguish between critical points lying in the same
orbit of the symmetry group $S_{(k)}$ (this group was defined in the {\it Step IV} of
Section~\ref{s23}).

Fix $[\kk]$ and divide the auxiliary variables $\T=\{t_l(j)\}$
(see Section~\ref{s23}) into  $(n-1)$  sets, in accordance with
the steps of iteration. Namely, the variables corresponding to the
first iteration are $\{t_i(j)\,, \  1\leq i\leq k_{1,\,j}\,, \
1\leq j\leq N\}$,  and to the $l$-th iteration are
$$
\{t_i(j)\,, \quad  k_{1,\,j}+\dots +k_{l-1,\,j}+1\,<\, i\, \leq k_{l+1,\,j}\,,\ \ 1\leq j\leq N\}\,,
\ \ 2\leq l\leq n-1\,.
$$
Perform the following change of variables. If  $t_i(j)$  belongs to the $l$-th set,
then we write
\begin{equation}\label{tau}
t_i(j)=s^{l+1}\tau_i(j)\,.
\end{equation}
If we re-write the function $\Psi$ in  variables $\tau=\tau_{[\kk]}=\{ \tau_i(j)\}$,
then the critical point system of
$\Psi(\tau)$ becomes a small deformation of the critical point system of the product
$$
\Psi_{[\kk]}\,=\,\Psi_1\cdot\, {\rm ...\,}\cdot\Psi_{n-1}\,,
$$
where every $\Psi_l$ is the master function corresponding
to $z=(0,1)$ and to singular vectors of the weight $\Lambda[\kk,\,l]$ in
$L_{\Lambda[\kk,\,l-1]}\,\otimes\, L_{\Lambda(l+1)}$\,; it depends on
variables $\{\tau_i(j)\}$ from the $l$-th set. Combining Theorem~\ref{T3}
and Proposition~\ref{PFu}, we conclude that the number of the (orbits of) critical
points with non-zero critical value of $\Psi_l$ is the maximal possible, i.e.
equals the dimension of the corresponding subspace of singular vectors,
for every $1\leq l\leq (n-1)$.

Any iterated Bethe vector $\vv_B$ corresponding to fixed $[\kk]$ is defined by
a certain critical point  of  the function  $\Psi_{[\kk]}$, and the $l$-th set of the
coordinates of this point is a critical point of $\Psi_l$, for every $2\leq i\leq n$. Hence,
according to the construction of  the basis of iterated Bethe vectors, the number
of the (orbit of)  critical points of $\Psi_{[\kk]}$ is exactly the number of the
iterated Bethe  vectors that correspond to $[\kk]$.
Notice, that near every critical point of  $\Psi_{[\kk]}$ there is at least one critical point
of the function  $\Psi( \tau)$.

These arguments work for every $[\kk]\in\KK$.
Taking into account that the whole number of  the (orbits of)  critical points of
the function $\Psi$ can not exceed  the dimension of the subspace of singular vectors
of the weight $\Lambda(\kk)$ in $L$, we conclude  that in fact a critical point
of function $\Psi(\tau_{[\kk]})$ lying near any critical point of  $\Psi_{[\kk]}$ is unique,
for every $[\kk]\in\KK$.

Denote this point  by $\tau^{(s)}$ and the  corresponding point of  $\Psi(\T)$ by
$\T^{(s)}$. If we substitute $\T^{(s)}$ into the universal weight function of
Section~\ref{s23}, then we obtain a Bethe vector $\vv^{(s)}=\vv(\T^{(s)})$.
If we make  change (\ref{tau})
and re-write $\vv(\T^{(s)})$  in terms of  $\tau^{(s)}$,  then all summands
coming from the iterated vector
$\vv_B$  will have  the same  factor $s^\nu$, where $\nu$ is a certain  integer,
whereas all other summands will include factors $s^a$ with $a<\nu$.     \hfill $\square$

\bigskip
The Lemma implies  that for $s$ big enough,  the  Bethe vectors
 in the Gaudin model associated with
 the weight $\Lambda(\kk)$ in $L$ and with $z=(s,s^2,\dots,\,s^n)$
 form a basis in the subspace of singular vectors of this weight in $L$.
The statement of the theorem then follows.   \hfill $\square$

\subsection{Basic Schubert intersections}\label{s42}
In order to calculate Bethe vectors, one should solve Bethe
systems. This is a difficult task, even for $n=2$; see~\cite{V} and \cite{MV2} for
a few of examples. Relations to the Schubert calculus suggest to 
look for non-degenerate planes in place of Bethe vectors.

\medskip
If $n=2$, then (\ref{I}) is the intersection of three Schubert varieties
and we call it a {\it basic Schubert intersection}.
By the same reason as in Section~\ref{s41}, one can always assume
$z_1=0$, $z_2=-1$.
We choose here $z_2=-1$ and not $z_2=1$  in order not to take care about
sings in  polynomials in our further calculation.

\medskip
\begin{prop}\label{PS} Let partitions $\w(0)$, $\w(-1)$ and $\w(\infty)$
satisfy
$$|\w(0)|+|\w(-1)|+|\w(\infty)|\,=\,\dim Gr_{N+1}(\Poly_d)\,.$$
If the intersection number
$\sigma_{\w(0)}\cdot\sigma_{\w(-1)}\cdot\sigma_{\w(\infty)}$ is positive,
then the basic Schubert intersection
$$
\Omega_{\w(0)}(0)\cap\Omega_{\w(-1)}(-1)\cap\Omega_{\w(\infty)}(\infty)
$$
is non-empty and the Wronskian of  any of it element is
$x^{|\w(0)|}(x+1)^{|\w(-1)|}$. \hfill $\square$
\end{prop}

\subsection{Special Schubert intersections}\label{s43}
Consider the basic Schubert intersections such that
at least one of  the Schubert varieties corresponding to $z_1,z_2$
is special. Without loss of generality we can and will
assume that this intersection has the form
\begin{equation}\label{ssi}
\I=\I_{\A,\,\w}=\Omega_\A(0)\cap\Omega_\w(\infty)\cap\Omega_{(m)}(-1)
\,\subset\, Gr_{N+1}(\Poly_d)\,,
\end{equation}
where $\A=(a_1,\dots\,,a_N)\,$ and $\w=(w_1,\dots\,,w_N)\,$ are
partitions, $d-N\geq a_1$, $d-N\geq w_1$, and $m=(N+1)(d-N)-|\A|-|\w|>0$.
We call such intersection a {\it special Schubert intersection}.

\medskip
If $V\in\I_{\A,\,\w}$, then the degrees
$d_1,\dots\,, d_N, d$ and the orders $0,\rho_1,\dots\,, \rho_N$ at $0$
of polynomials in $V$ are related to $\w$ and $\A$ as follows, see (\ref{rho}),
\begin{equation}\label{dr}
d_l=d-w_l+l-N-1\,,\ \ \rho_{N+1-l}=a_l+N+1-l\,,\ \ 1\leq l\leq N.
\end{equation}

The  famous Pieri formula (\cite{Fu}) can be reformulated now in terms of degrees
and orders at $0$ of  polynomials.

\begin{lemma}\label{PF} {\rm (Pieri formula)} If a
special Schubert intersection  {\rm (\ref{ssi})}
is non-empty, then it consists of a single element $V$, and
\begin{equation}\label{drho}
0\leq d_1<\rho_N\leq d_2<\rho_{N-1}\leq d_3<\,\dots\,<\rho_1\leq d\,,
\end{equation}
where $d_1,\dots\,, d_N, d$ are the degrees and $0,\rho_1,\dots\,, \rho_N$
are the orders at $0$ of the polynomials in $V$. The Wronskian of $V$ is
$$
W_V(x)=x^{|\A|}(x+1)^{m}=x^{\rho_1+\dots+\rho_N-N(N+1)/2}(x+1)^m\,.
$$
\ \hfill $\square$
\end{lemma}

The relation (\ref{drho}) says that the polynomial of degree $d_1$ in $V$
has order $0$ at $0$. Hence there is a polynomial of degree $d_2$ in $V$
that has order $\rho_N$ at $0$, and so on. We arrive at the following conclusion.

\begin{cor}\label{cQ}
 If a special Schubert intersection  {\rm (\ref{ssi})} is non-empty, then
it is $\Span \{Q_1(x)\,,\dots,\,Q_{N+1}(x)\}$,  where $Q_i(x)$
has degree $d_i$  and  order $\rho_{N+2-i}$   at $0$ for every
$1\leq i\leq N+1$. Here $d_i$ and $\rho_i$  are defined for $1\leq i\leq N$
by \rm{(\ref{dr})},\,  $d_{N+1}=d$, $\rho_{N+1}=0$.  \hfill $\square$
\end{cor}

\section{A basis of the special Schubert intersection}\label{s5}

In this section we produce explicitly a basis of the (single element of the)
special Schubert intersection $V=\I_{\A,\,\w}$ given by (\ref{ssi})--(\ref{drho}).
Section~\ref{s51} is devoted to the simplest case, when the value of $m$ is the maximal
possible, i.e. $m=d-N$; in Section~\ref{s52} the generic case is done.

\subsection{Truncated binomials }\label{s51}
Let $m,d$ be positive integers, $m\leq d$. We call the polynomial
$$
P_{m;\,d}(x):=1+dx+{d\choose 2}x^2+\dots +{d\choose m}x^m
$$
a {\it truncated binomial}. In particular, $P_{d;\,d}(x)=(x+1)^d$.

\begin{thm}\label{T4} The truncated binomials
$$P_{m_1;\,d}(x)\,,\ \dots\,,\ P_{m_N;\,d}(x)\,,\ P_{d;\,d}(x)\,,
$$
where $0\leq m_1<m_2<\dots<m_N<d$, span the (single element of the)
Schubert intersection
$$
\I_{\A,\,\w}\,=\,\Omega_{\A}(0)\cap\Omega_{(d-N)}(-1)\cap
\Omega_{\w}(\infty)\,\subset\,Gr_{N+1}(\Poly_d)\,,
$$
where
$$
\begin{array}{l}
\A\,=\,(m_N+1-N, \dots\,, m_2-1,m_1 )\,,\\
\w \,=\,(d-m_1-N, \dots\,,d-m_{N-1}-2, d-m_N-1)\,.
\end{array}
$$
\end{thm}

\medskip\noindent{\bf Proof.}\ \ Notice that
$|\A|+|\w|+d-N=\dim Gr_{N+1}(\Poly_d)$,
therefore $\I_{\A,\,\w}$ consists of at  most one
element (by the  Pieri formula). Thus it is enough to check that
$$V=\Span\{P_{m_1;\,d}(x)\,,\ \dots\,,\ P_{m_N;\,d}(x)\,,\ P_{d;\,d}(x)\}$$
belongs to this Schubert intersection.
Consider polynomials
$$
Q_1(x)=P_{m_1;\,d}(x)\,, \ Q_i(x)=P_{m_i;\,d}(x)-P_{m_{i-1};\,d}(x)\,,
\ \ 2\leq i\leq N+1\,.
$$
These polynomials also span $V$. Their orders  at $0$ are $0,
m_1+1,..., m_N+1$, and their orders at $\infty$ are $0, d-m_N,...,
d-m_1$  (compare with Corollary~\ref{cQ}). Moreover the sum
$$Q_1(x)+\dots+Q_{N+1}(x)\,=\,P_{d;\,d}(x)\,\in \,V
$$
has order $d$ at $(-1)$. Hence the other orders of polynomials of $V$ at $(-1)$
have to be $0,1,...,N-2$. Indeed, otherwise the intersection of the
corresponding Schubert varieties would have a negative dimension.
Thus $V$ belongs to  the three  Schubert varieties. \hfill $\square$

\begin{cor}\label{cP}\ \
If $0\leq m_1<\dots\,<m_N<d$, then the Wronski determinant of
truncated binomials $P_{m_1;\,d}(x),\dots\,,P_{m_N;\,d}(x),P_{d;\,d}(x)$
is
$$
\det\left(\frac{d^j}{dx^j}P_{m_i;\,d}(x)\right)= c\cdot
x^{m_1+\dots\,+m_N-N(N-1)/2}(x+1)^{d-N}\,,\ \ 0\leq i,j\leq N\,,
$$
where $c$ is a non-zero constant.  \hfill $\square$
\end{cor}

\medskip
  If $\Span\{P_{m_1;\,d},\dots\,,\ P_{m_N;\,d},\ P_{d;\,d}\}$
 is non-degenerate, then the truncated binomials define a Bethe vector
 of the weight
 $$\w^*\,=\,(d-m_1-N)\lambda_1\,+\,\sum_{i=2}^N\,(m_i-m_1-i+1)\lambda_i
 $$
 in the Gaudin model associated with $(z_1,z_2)=(0,-1)$ and
 the tensor product $L_{\Lambda(0)}\otimes L_{\Lambda(-1)}$ of
irreducible $sl_{N+1}$-representations with highest  weights
 $$
 \Lambda(0)\,=\,\sum_{i=1}^N\,(m_{N+1-i}-N+i)\lambda_i\quad {\rm and}\quad
 \Lambda(-1)\,=\,(d-N)\lambda_1\,,
 $$
 respectively, see Corollary~\ref{C1}.

\medskip
For $N=1$ the plane $\Span\{P_{k;\,d},\ P_{d;\,d}\}$ is non-degenerate, for any $k<d$,
see~\cite{Sc2} . In fact this follows from an elementary property of the truncated binomials.

\begin{lemma}\label{lP}
 $P_{k;\,d}(x)$ does not have multiple roots for $ k<d$.
\end{lemma}

\medskip\noindent{\bf Proof.}\ \ We have  $P_{k;\,d}(0)=1.$ On the other hand,
$$
d\cdot\,P_{k;\,d}(x)\,-\,(1+x)P'_{k;\,d}(x)\,=\, d\cdot\, {d-1\choose k}x^k\,;
$$
hence if  there exists a common root of $P_{k;\,d}$ and $P'_{k;\,d}$\,, then it
is $x=0$.   \hfill $\square$

\medskip
For $N>1$, we do  not know if  $\Span\{P_{m_1;\,d},\dots\,,\ P_{m_N;\,d},\ P_{d;\,d}\}$
is always non-degenerate.  Notice that for $N=2$ the only additional requirement is as follows.

\begin{conj}\label{CP} Truncated binomials
 $P_{k;\,d}(x)$ and $P_{m;\,d}(x)$ do not have
common roots for $0< k<m< d$\,.
\end{conj}

\medskip
Indeed, for $V=\Span\{P_{k;\,d}(x)\,,\ P_{m;\,d}(x)\,,\ P_{d;\,d}(x)\}\,$ the  flag
$V_{\bullet}=\{V_1\subset\,V_2\subset\,V\}$ defined by (\ref{V})  in
Section~\ref{s33}  is as follows,
$$
V_1=\Span\{P_{k;\,d}(x)\}\,,\quad V_2= \Span\{P_{k;\,d}(x)\,,\
P_{m;\,d}(x)\}\,,
$$
and the intermediate Wronskians are
$$
W_1(x)=W_{V_1}(x)=P_{k;\,d}(x)\,,\quad
W_2(x)=W_{V_2}(x)=P_{k;\,d}'(x)P_{m;\,d}(x)-P_{k;\,d}(x)P_{m;\,d}'(x)\,.
$$
By definition (see Section~\ref{s34}), $V$ is non-degenerate if
\begin{itemize}
\item
$P_{k;\,d}(x)=W_1(x)$ do not have multiple roots;
\item
$W_1(x)$ and $W_2(x)$ do not have
common roots distinct from $0$ and $(-1)$.
\end{itemize}
The first item holds (it is Lemma~\ref{lP}), therefore the second item
implies that every common root of $W_1(x)$ and $W_2(x)$ should be
a common root of $P_{k;\,d}(x)$ and $P_{m;\,d}(x)$ as well.

\medskip
The statement of the conjecture looks rather elementary, and holds in
examples (see  Section~\ref{s63}).
However we did not succeed either prove or disprove it.

\subsection{A basis in $\I_{\A,\,\w}$}\label{s52}
Consider a special  Schubert intersection that satisfies
(\ref{ssi})--(\ref{drho}). It is either
empty  or consists of a single plane $\I_{\A,\,\w}= V\in Gr_{N+1}(\Poly_d)\,,$
according to Lemma~\ref{PF}.
Let us begin with the relation $V\in \Omega_{(m)}(-1)$. It means
that $V$ contains a polynomial of order $m+N$ at $(-1)$, 
\begin{equation}\label{P}
P(x)=(x+1)^{m+N}+C_1(x+1)^{m+N+1}+\,\dots\,+C_{d-m-N}(x+1)^d\,\in V\,,
\end{equation}
where  $C_1,\dots\,, C_{d-m-N}$ are some  constants.

\begin{prop}\label{det}\ \
Constants $C_1,\dots\,, C_{d-m-N}$ are uniquely defined by the special
Schubert  intersection  $\I_{\A,\,\w}$ given by \rm{(\ref{ssi})--(\ref{drho})}.
\end{prop}

\medskip\noindent{\bf Proof.}\ \  The relation (\ref{drho}) of  Lemma~\ref{PF}
says that $P(x)$ does not contain terms $x^l$ with $l\in\LL$, where
\begin{equation}\label{l}
\LL\,=\,\{\,d_1+1,\dots\,,\,\rho_N-1,\,d_2+1,\dots\,,\,\rho_{N-1}-1,
\dots\,,\,d_N,\dots\,,\,\rho_1-1\,\}.
\end{equation}
Denote by $l_i$  the $i$-th entry in $\LL$. We get
 $d-m-N$ linear equations on  $C_1\,,\dots\,, C_{d-m-N}\,$,
\begin{equation}\label{C}
{m+N \choose l_i}\,+\,C_1\cdot {m+N+1 \choose l_i}\,+\,\dots\,+\,
C_{d-m-N}\cdot {d \choose l_i}\,=0\, ,\quad  1\leq i \leq d-m-N\, .
\end{equation}
The determinant of this system coincides, up to a sign, with
$$
A\,=\,\det\,\left(A_{ij}\right)\,,
$$
where  $A_{ij}$'s are binomial coefficients,
$$
A_{ij}\,=\,{d+1-j \choose l_i}\,,\ \  1\leq i,j \leq d-m-N=K\,.
$$

First of all, we see that the $i$-th row has $\left(l_i!\,(d+1-K-l_i)!\right)^{-1}$
 as a common factor, and the $j$-th column   $(d+1-j)!$, i.e.
$$
A\,=\,\frac{\prod_{j=1}^K\, (d+1-j)!}
{\prod_{i=1}^K\,\left(l_i!\,(d+1-K-l_i)!\right) }\, \cdot\, B\,,
$$
where
$$
B\,=\,\det\,\left(B_{ij}\right)\,; \quad
B_{ij}\,=\,\frac1{(d-l_i-K+2)\,\dots\, (d-l_i+1-j)}\ \ \ {\rm for}
\   j\neq K;  \ \ \  B_{iK}=1\,.
$$
Introduce notation $b_i=d-l_i-K+2$, for $1\leq i\leq K$. With this notation,
we have
$$
B_{ij}\,=\,\frac1{b_i\cdot\,(b_i+1)\dots\,(b_i+K-1-j)} \ \ \ \ {\rm for}
\   j\neq K\,.
$$
Next, we apply  to $B_{ij}$ the following simple fraction decomposition.

\begin{lemma}\ \ For any positive integer $M$, we have
$$
\frac1{x\cdot\,(x+1)\dots\,(x+M)}\,=\,\sum_{i=0}^M\,
\frac{(-1)^i}{i!\,(M-i)!\, (x+i)}\,.
$$
\ \ \   \hfill  $\square$
\end{lemma}

Thus using the $(K-1)$-th column, which consists of $1/b_i$, we can erase the summands $1/b_i$ in the all previous $K-2$ columns. Then the $(K-2)$-th column becomes $B_{i,K-2}=-1/(b_i+1)$, and we can use it to erase $1/(b_i+1)$ in the
all previous $K-3$ columns and so on. We get
$$
B=\frac{(-1)^{(K-1)(K-2)/2}}{1!\,2!\,\dots\, (K-2)!}\cdot \tilde B\,,
$$
where
$$
\tilde B\,=\,\det\,\left(\tilde B_{ij}\right)\,,\quad \tilde B_{ij}\,=\,
\frac1{b_i+K-1-j}\ \ {\rm for}\   j\neq K;  \ \ \  \tilde B_{iK}=1\,.
$$
To calculate this determinant is an exercise in linear algebra.
Write $\tilde B\,=\,\tilde B_K(b_1,\dots,\,b_K)$.

\begin{lemma}\ \ We have the following recursion,
$$
\tilde B_K(b_1,\dots,\,b_K)\,=\,\frac{(-1)^{(K-1)(K-2)/2}\,(K-2)!\,
\prod_{i=1}^{K-1}\,(b_k-b_i)}{\prod_{i=1}^K\,b_i\,
\prod_{j=1}^{K-2}\,(b_K+j)}\,\cdot\, \tilde B_{K-1}\,(b_1+1,\dots,\, b_{K-1}+1)\,.
$$
\ \ \   \hfill  $\square$
\end{lemma}

Finally we obtain
$$
\det\,\left(A_{ij}\right)=\pm\, \prod_{i=1}^{K-1}\, (d+1-i)^i\, \cdot \,
\prod_{i=1}^K\, {d+1-K \choose l_i}\, \cdot \,
\prod_{1\leq j<i\leq K}\,(l_j-l_i)\,,
$$
and hence  $\det\,\left(A_{ij}\right)\,\neq\, 0$.  \hfill $\square$

\bigskip
Now  substitute the solution $C_1,\dots\,, C_{d-m-N}$ of the system (\ref{C})
into (\ref{P}).  Recall that we set $d_{N+1}=d$ and $\rho_{N+1}=0$.
The polynomial  $P(x)$ can be presented as the sum of $N+1$ polynomials,
$$
P(x)\,=\,Q_1(x)+\,\dots\,+Q_{N+1}(x)\,,
$$
such that for every $1\leq i\leq N+1$ the polynomial $Q_i(x)$ has degree
$d_i$ and order $\rho_{N+2-i}$ at $0$.
This means that
$$
\Span\{Q_1(x),\,\dots\,, \,Q_{N+1}(x)\}\,\in\,
\Omega_{\A}(0)\cap\Omega_{(m)}(-1)\cap\Omega_{\w}(\infty)
$$
 (compare with Corollary~\ref{cQ}). We arrive at the following result.

\begin{thm}\label{T5} Every special  Schubert intersection $\I_{\A,\,\w}$
given by {\rm (\ref{ssi})--(\ref{drho})} is non-empty, and
the polynomials   $Q_1(x),\,\dots\,, \,Q_{N+1}(x)$ defined above provide a
basis of (the single element of) $\I_{\A,\,\w}$. \hfill $\square$
\end{thm}

\subsection{Examples: bases of special Schubert intersections and
Bethe vectors}\label{s63} If a basis of  an element in
$Gr_{N+1}(\Poly_d)$ is known, then one can find the additional
roots of the intermediate  Wronskians and to check the
non-degeneracy conditions, as it is explained in
Section~\ref{s34}. If these conditions hold, then the substitution
of  the additional roots into the universal weight function
$\vv(\T)$ ( see Section~\ref{s23}) gives the corresponding Bethe
vector. All special Schubert intersections in our examples are
non-degenerate.

\medskip
First we discuss special Schubert intersections given by presentations
of $(x+1)^5$ as a sum of polynomials.

\bigskip\noindent
{\bf 1.}\ \ Polynomials
$$
Q_1(x)=1+5x\,,\ \ Q_2(x)=10x^2+10x^3\,,\ Q_3(x)=5x^4+x^5
$$
span the special Schubert intersection $\I_{\A,\w}\in Gr_3(\Poly_5)$, where
$\A=\w=(2,1)$ (as Theorem~\ref{T4} or an easy direct calculation shows).
The corresponding Bethe vector is a highest
weight vector of the weight $\Lambda_{\w^*}=2\lambda_1+\lambda_2$
in the tensor product $L_{\Lambda(0)}\otimes\,L_{\Lambda(-1)}$ of
$sl_3$-representations, where $\Lambda(0)=\Lambda(-1)=2\lambda_1+\lambda_2$.
We have
$$
\Lambda_{\w^*}\,=\,\Lambda(0)+\Lambda(-1)-2\alpha_1-\alpha_2.
$$
In order to find the Bethe vector, we have to calculate values of the
variables $t_1(1)$, $t_2(1)$ corresponding to $\alpha_1$ and  of
the variable $t(2)$
corresponding to $\alpha_2$, according to Corollary~\ref{C1}.
These values are the additional roots of the Wronskians
$W_2(x)=Q_1(x)Q'_2(x)-Q'_1(x)Q_2(x)$ and $W_1(x)=Q_1(x)$, respectively,
as Corollary~\ref{C1} says. We have (up to a constant factor)
$$
W_1(x)=1+5x\,,\ \ W_2(x)=10x(10x^2+8x+2)\,,
$$
i.e. the value of $t(2)$ is $-1/5$ and the values of $t_{1,2}(1)$ are $(-2\pm i)/5$.
The substitution into the function $\vv(\T)$ defined in Section~\ref{s23}
gives the Bethe vector.

\bigskip\noindent
{\bf 2.}\ \ Another presentation of $(x+1)^5$ as the sum of
three polynomials,
$$
Q_1(x)=1\,,\ \ Q_2(x)=5x+10x^2+10x^3+5x^4\,,\ Q_3(x)=x^5\,,
$$
gives $\Span\{Q_1,Q_2,Q_3\}=\I_{\A,\w}\in Gr_3(\Poly_5)$ with
$\A=\w=(3,0)$. Now $\w^*=(3,3)$ and the corresponding Bethe vector
has the weight
$$
\Lambda_{\w^*}=3\lambda_1+3\lambda_2=6\lambda_1-3\alpha_1=
\Lambda(0)+\Lambda(-1)-3\alpha_1\,.
$$
There are no variables corresponding to $\alpha_2$ (indeed, $W_1(x)=Q_1(x)=1$
has no roots), and there are three variables corresponding to $\alpha_1$.
The values of these variables that give the Bethe vector are the roots of
$W_2(x)=Q'_2(x)=5(4x^3+6x^2+4x+1)$.

\bigskip\noindent
{\bf 3.}\ \ If we take
$$
Q_1(x)=1\,,\ \ Q_2(x)=5x+10x^2\,,\ \ Q_3(x)=10x^3+5x^4\,,\ Q_4(x)=x^5\,,
$$
we get $\Span\{Q_1,Q_2,Q_3,Q_4\}=\I_{\A,\w}\in Gr_4(\Poly_5)$ with
$\A=(2,1,0)$ and $\w=(2,1,0)$, $\w^*=(2,2,1)$. We have
$$
\Lambda_{\w^*}\,=\,\Lambda(0)+\Lambda(-1)-2\alpha_1-\alpha_2.
$$
The corresponding Bethe vector is $\vv(\T^0)$ where $t^0(2)$
is the root of $W_2(x)=Q'_2(x)$ and $t^0_{1,2}(1)$ are
the non-zero roots of $W_3(x)=Q'_2(x)Q''_3(x)-Q''_2(x)Q'_3(x)$.

\bigskip
Now consider  $(x+1)^4+c(x+1)^5$ and choose $c$ in such
a way that $x^3$ disappears. An easy calculation gives $c=-2/5$.

\bigskip\noindent
{\bf 4.}\ \ Present the obtained polynomial as the  following  sum,
$$(x+1)^4-\frac 25(x+1)^5=\left(\frac 35 +2x\right)\, +\, 2x^2 -
\left(x^4+\frac 25x^5\right)\,=\,Q_1(x)+Q_2(x)-Q_3(x)\,.
$$
According to Theorem~\ref{T5}, we have
$\Span\{Q_1,Q_2,Q_3\}=\I_{\A,\w}\in Gr_3(\Poly_5)$
where $\A=(2,1)$, $\w=(2,2)$ (of course it is easy to check this directly).
Hence $\w^*=(2,0)$ and in order to obtain the Bethe vector of the
weight $2\lambda_1=\Lambda(0)+\Lambda(-1)-\alpha_1-\alpha_2$ one has to
substitute the roots of two linear equations,
$3/5+2x=0$ and $3/5+x=0$, into the weight function $\vv(\T)$ of Section~\ref{s32}.

\begin{remark}\label{last} {\rm
In more complicated examples, with larger numbers $N$ and $d$,
it is still easy to calculate the elementary symmetric functions in the
auxiliary variables $t(i)$, that is the coefficients of the polynomials
$T_{N+1-i}(x)$. Notice that the functions $\omega(\T)$ entering the universal
weight function, see (\ref{omega}), are in fact functions in the same
elementary symmetric functions, as  Step IV of Section~\ref{s23} shows.
It would be helpful to re-write $\vv(\T)$ in terms of polynomials
$T_{N+1-i}(x)$. With S.~Chmutov, we obtained the corresponding expressions
for some of $\omega(\T)$'s in \cite{CSc}.}
\end{remark}

\end{document}